\documentclass[12pt]{article}
\usepackage{epsfig,graphicx,graphics}
\usepackage{epsfig, wrapfig} 
\usepackage{rotating}
\newcommand{\fracd}[2]{  \frac {\displaystyle {#1}}{\displaystyle {#2}}   }
\newcommand{\bb}{\begin{equation}}
\newcommand{\ee}{\end{equation}}

\newcommand{\esp}{\vspace*{0.5 cm}}

\begin{document}

\begin{center}
{\Large {\bf A MATHEMATICAL APPROACH TO THE PLATO'S PROBLEM }}
\end{center}

\esp

\begin{center}
{\bf Luiz Bevil\'acqua, Adilson J. V. Brand\~ao and  Rodney C. Bassanezi

Universidade Federal do ABC - UFABC, Santo Andr\'e - SP, Brazil \footnote{luiz.bevilacqua@ufabc.edu.br, adilson.brandao@ufabc.edu.br, rodney.bassanezi@ufabc.edu.br}
}
\end{center}

\esp

{\bf ABSTRACT.} Maybe the first inverse problem presented in the history of the occidental thought is described in the book {\it Republic}, written by Plato. The problem is posed in the Book VII in a text known as the {\it Allegory of the Cave}. That text motivated us to formulate a simple mathematical model that simulates, in a sense, the situation of the persons described in that problem.

\esp

\renewcommand{\theequation}{1.\arabic{equation}}
\setcounter{equation}{0}

{\bf 1. INTRODUCTION. } The motivation of this paper is a problem posed more than 2000 years ago by Plato. Nevertheless the essence of that problem is one of the most challenging questions to be found in modern science and technology. Plato in the Book VII of his monumental masterpiece {\it Republic}, wrote a dialogue known as the {\it Allegory of the Cave} \cite{Plato}. That text is maybe the first inverse problem presented in the history of the occidental thought. Briefly, several persons, sitting, facing the rear wall of a cave were totally immobilized even unable to move their heads. A fire was set between the persons and the wall. A procession of men and women of all ages and animals carrying, pulling and pushing all kind of objects moved between the wall and the fire. The persons tied on the chairs were able only to watch the images projected on the wall but could not see the parade directly. The challenge was to interpret the images and find out what kind of persons, animals and objects composed the parade. Plato explored through this allegory the knowing process and how to reach the truth. Despite the fact that his aim was distinct he presented maybe for the first time in a well organized structure an inverse problem, specifically, pattern recognition. The {\it Allegory of the Cave} is one of the most marvelous texts ever produced by the human mind.       

Plato's text motivated the design of the kit presented in the sequel that simulates in a simplified version the images observed by the persons tied to the chairs in Plato's story. The first part of this paper will describe the direct problem and the last section will deal with a particular inverse problem. 

Nowadays one of the most important challenges that falls into this category is the pattern recognition problem. There are several applications in medical diagnosis through the interpretation of images obtained with sophisticated equipments as computerized tomography and nuclear magnetic resonance devices. The question is: given the pattern what is the tissue density and distribution in the organ under examination. 

The problem posed here is essentially the same except that the number of parameters involved is substantially smaller.  
If we know how to solve the direct problem the solution of the corresponding inverse problem can be substantially facilitated.  For the reader interested in more information on inverse problems see \cite{Tarantola}.

\renewcommand{\theequation}{2.\arabic{equation}}
\setcounter{equation}{0}

{\bf 2. THE PROBLEM. } A slender object with the shape of a thin plate $B$ rotates with angular velocity $\Omega$ around an axis $O$ perpendicular to a plane containing the object. Such an object is confined inside a box $C$ such that the slot $F$ and the object belong to the same plane.
\begin{figure}[!ht]
\centering
\includegraphics[scale=0.5]{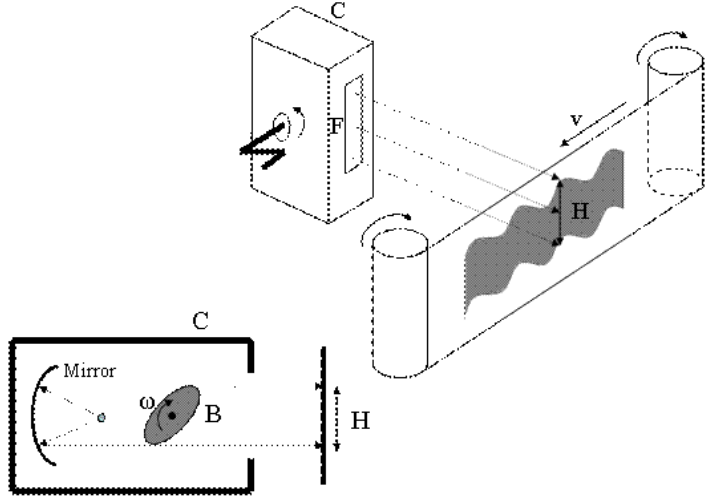}
\caption{Building setup: object and projection system}
\end{figure}

A device inside the box, consists of a lamp and a concave mirror which projects the profile of the plate $B$ on a moving screen $E$. This projection produces a shadow with height $H$  equal to the distance between the upper and lower tangents to the contour of the plate $B$.

If the film is sensitive to light, and unwinds in front of the slot with velocity $v$ at the same time that the plate rotates around $O$, an image which is in a certain way a "picture of $B$" will be generated.
\begin{figure}[!ht]
\centering
\includegraphics[scale=0.5]{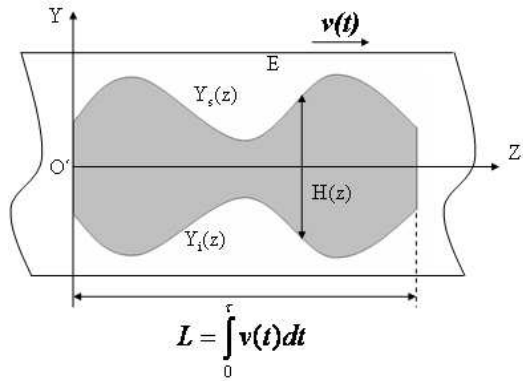}
\caption{Image of the object B recorded on the film}
\end{figure}

The image printed on the tape is a continuous picture with the shape of a ribbon with variable width, we will call it "kinematic image" of $B$.

The picture is bounded above and below by $Y_s (z)$ and $Y_i (z)$, respectively. The width is given by 
$H(z)=Y_s (z) - Y_i (z)$.

The following problems arise:
\begin{enumerate}
\item[i)] Given the geometry of $B$, the angular velocity $\Omega (t)$ and the travelling velocity $v(t)$ find $H(z,t_0)=H_0 (z)$ after a given time $t_0$.
\item[ii)] Given $Y_s (z), Y_i (z)$ and $H_0 (z)$ find the geometry of $B$. Is the solution unique? 
\end{enumerate}

The first problem is known as a "direct problem" and the second one is an "inverse problem".

{\bf Observations:} 
\begin{enumerate}
\item The projection $H_0 (z)$ do not identify reentrancies. In other words, it is a critical condition the contour to be convex in order to have an unique solution.
\begin{figure}[!ht]
\centering
\includegraphics[scale=0.5]{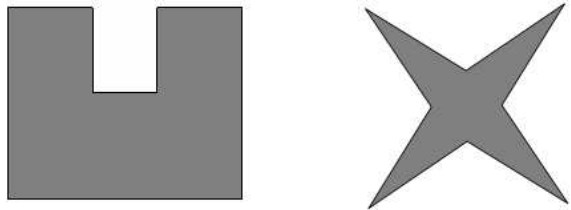}
\caption{Non-convex contours}
\end{figure}

\item Another difficulty is the existence of vertices on the contour. These singularities could be overcome, but they will need special treatment since at the vertices there is not a definite tangent. We will see that for the case of regular polygons it is possible to determine a simple solution.
\begin{figure}[!ht]
\centering
\includegraphics[scale=0.5]{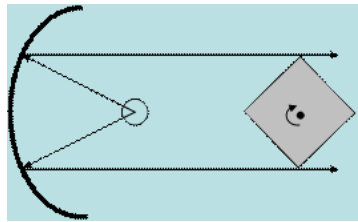}
\caption{Contour with vertices: tangents not well-defined}
\end{figure}

\end{enumerate}

\renewcommand{\theequation}{3.\arabic{equation}}
\setcounter{equation}{0}

{\bf 3. FORMULATION OF THE DIRECT PROBLEM.} Consider a flat object $B$ with a smooth contour. Let the contour, be defined in terms of polar coordinates $(\beta,r(\beta))$ having as a reference $Qxy$ fixed in $B$. The object $B$ moves with respect to $OXY$ rotating around a perpendicular axis to its plane, going through the point $O$ fixed in $OXY$, with angular velocity $\theta'(t)=d\theta /dt$, where $\theta (t)$ is the angle between a line $Qx$ fixed in $B$ and a line $OX$ fixed at the inertial reference system.
\begin{figure}[!ht]
\centering
\includegraphics[scale=0.5]{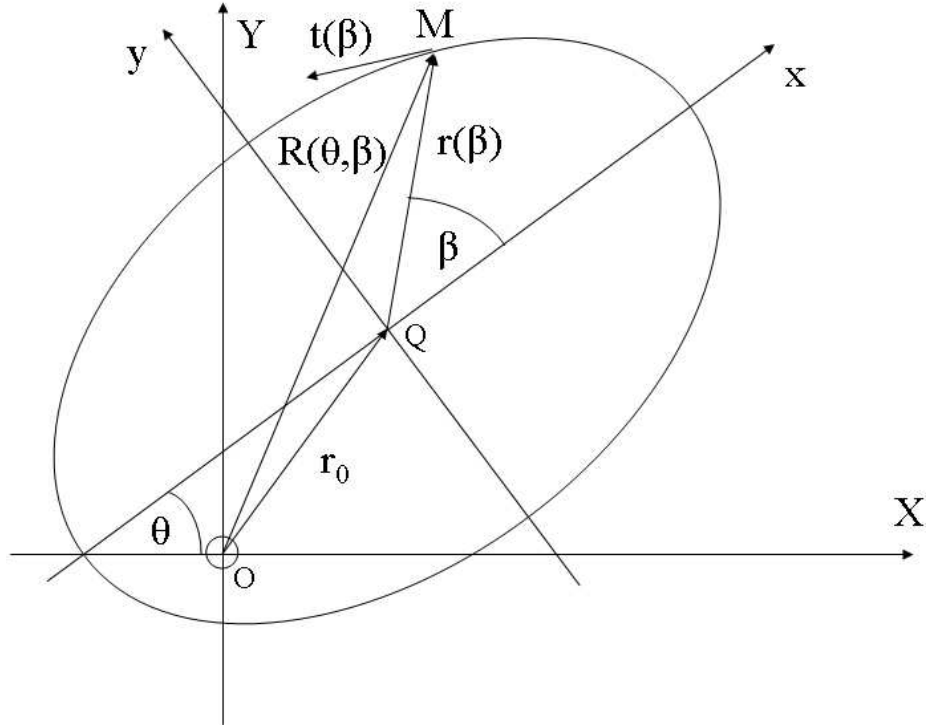}
\caption{Contour and tangent representation}
\end{figure}

Our propose is to find the shadow of $B$ projected on the rotating tape, that is, the distance $H(z,t)$ between the upper and lower tangents.
\begin{figure}[!ht]
\centering
\includegraphics[scale=0.5]{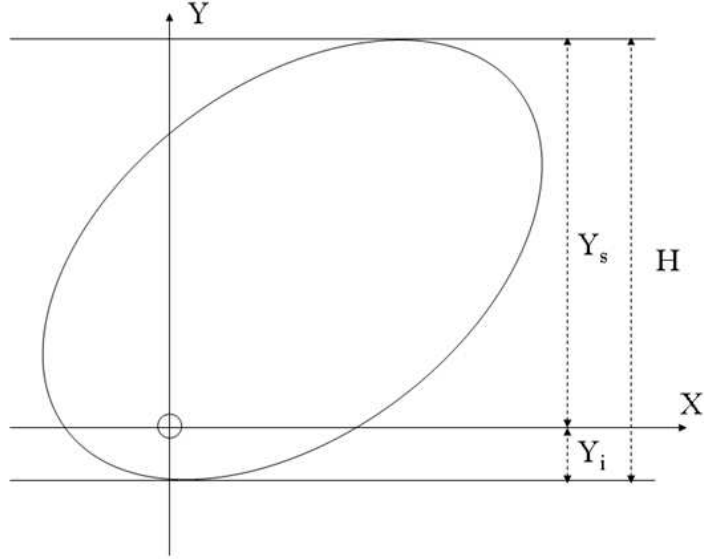}
\caption{Limits $Y_s$ and $Y_i$ of the projection and height $H$}
\end{figure}

Let $\vec{r}(\beta)$ be the position vector of an arbitrary point $M$ at the contour of $B$. We are admiting a smooth contour, that is, without corners. The general expression of the tangent related to the contour of $B$ is:
\[
\vec{t}= \fracd{d\vec{r}}{d\beta}.
\]

Clearly the upper and lower tangents are paralell to the horizontal axis $OX$. This is equivalent to say that both are orthogonal do $OY$. Therefore we may write:

\begin{equation}\label{eq1}
\vec{t}\cdot \vec{j}=0,
\end{equation}

where $\vec{j}$ is the unitary vector corresponding to the axis $OY$.

The representation of $\vec{t}$ in $OXY$ is given by $\vec{T}={\cal R} \vec{t}$ where the rotation matrix is:

\[ {\cal R}=\left [
\begin{array}{cc}
	\cos \theta & -\sin \theta \\
	\sin \theta & \cos \theta \\
\end{array}
\right ].
\]

The condition (\ref{eq1}) can therefore be written as:

\begin{equation} \label{eq2}
{\cal R} \vec{t} \cdot \vec{j}=0.
\end{equation}

If the contour is smooth and convex there will be two solutions for this equation corresponding to the upper and lower tangents. This expression relates $\beta$ and $\theta$ in such a way that with this relationship $\beta$ will always define the same position of the extreme tangents. We will designate this position $\beta (\theta) = \beta^*$.

In general there is not an explicit form $\beta (\theta)$, and hence, $\beta^*$ has to be determined from the numerical solutions of (\ref{eq2}). In the next sections we will to study some special cases where it is possible to obtain $\beta^*$ explicitly. 

Once we determine $\beta^*$, or $\beta_u^*$ and $\beta_l^*$ corresponding to the upper and lower tangents it is possible to find $Y_s$ and $Y_i$:

\[ 
\begin{array}{c}
Y_s = \vec{R}(\beta_u^*,\theta) \cdot \vec{j}=\vec{R}_u \cdot \vec{j}, \\
Y_i = \vec{R}(\beta_l^*,\theta) \cdot \vec{j}=\vec{R}_j \cdot \vec{j}, 
\end{array}
\]

where $\vec{R}_u$ and $\vec{R}_l$ are position vectors defining the tangent points $M$ and $N$ on the contour corresponding to the upper and lower tangents represented in $OXY$, that is:

\[ 
\begin{array}{c}
\vec{R}_u = {\cal R} (\vec{r}_0 + \vec{r}(\beta_u^*)), \\
\vec{R}_l = {\cal R} (\vec{r}_0 + \vec{r}(\beta_l^*)). 	
\end{array}
\]

Now, 
\[
H=Y_s - Y_i = [{\cal R}(\vec{r}(\beta_s^*) - \vec{r}(\beta_s^*))]\cdot \vec{j}, 
\]

that is, the total height of the shadow does not depend upon the pole of rotation $O$.

{\bf Example:} The figure 7 shows two projections of the same circle with different poles.
\begin{figure}[h]
\centering
\includegraphics[scale=0.5]{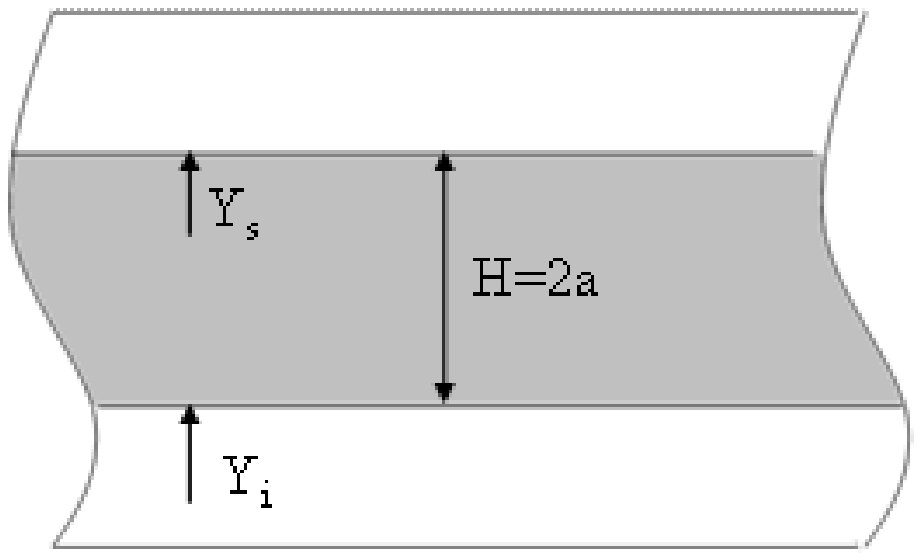}
\includegraphics[scale=0.5]{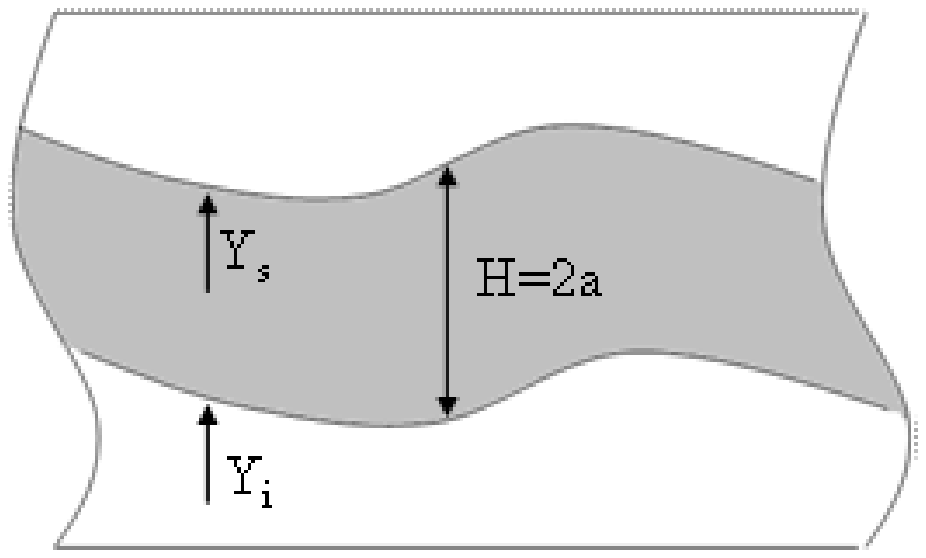}
\caption{Image of a circle. Rotation around the center and around a point on the contour}
\end{figure}

It is now necessary to find the image printed on the moving film that travels in front of the small hole with velocity $v(t)$.

Consider the axis $O'Z$ on the film, orthogonal to the $OXY$ plane. Let $O'$ be the origin of the $Z$ axis. Clearly the figure printed on the film, that is, the kinematic image of $B$, will be comprised between $Y_s$ and $Y_i$, corresponding to the points where the upper and lower tangents intercept the tape.

The projection $Y_s$ on the film varies with time as:

\[
\frac{dY_s}{dt}=\frac{d}{dt}[{\cal R} (\vec{r}_0 + \vec{r}(\beta_u^*))\cdot \vec{j}]=\frac{\partial}{\partial \theta}
[{\cal R} (\vec{r}_0 + \vec{r}(\beta_u^*))\cdot \vec{j}]\frac{d\theta}{dt}.
\]

On the other hand the film moves with speed $v$, therefore
\[
\frac{dz}{dt}=v(t).
\]
The solution of the system
\[ 
\begin{array}{cl}
\fracd{dY_s}{dt}=& \frac{d}{dt}[{\cal R} (\vec{r}_0 + \vec{r}(\beta_u^*))\cdot \vec{j}],\\
\fracd{dz}{dt} =& v(t) 	
\end{array}
\]
gives $Y_s = Y_s (z)$. Similarly we would obtain $Y_i = Y_i (z)$.

\begin{figure}[h]
\centering
\includegraphics[scale=0.4,viewport=130 80 589
482,clip]{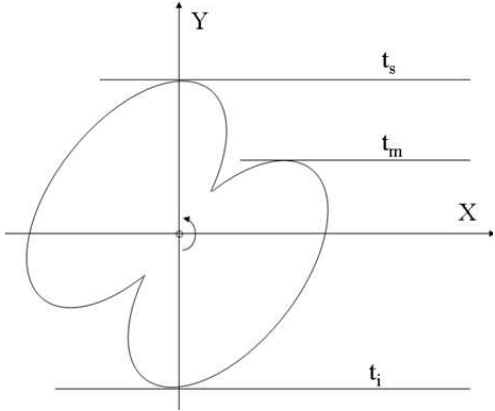}
\caption{Non-convex contour with more than two tangents}
\end{figure}

Note that the solution of the equation (\ref{eq2}) has two non-coincident roots if and only if the object $B$ is convex and its contour is a smooth curve.

\newpage

\renewcommand{\theequation}{4.\arabic{equation}}
\setcounter{equation}{0}

{\bf 4. THE CASE OF THE CIRCLE.}
Let a circle of radius $a$ rotates around the centre $O$.
\begin{figure}[!ht]
\centering
\includegraphics[scale=0.5]{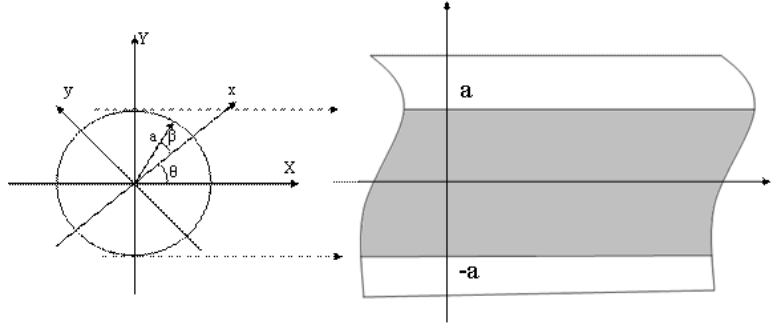}
\caption{Image of a circle rotating around the center}
\end{figure}

The solution is trivial. The kinematic image is a strip with width equal $2a$, bounded by  $Y_s = a$ and $Y_i=-a$. Formally:
\[ 
\vec{r}= a \left [
\begin{array}{cc}
 \cos \beta \\
 \sin \beta  \
\end{array} \right ],
\vec{t}= \left [
\begin{array}{cc}
 -\sin \beta \\
 \cos \beta  \
\end{array} \right ],
\]

\[{\cal R} \vec{t}=\left [
\begin{array}{cc}
	\cos \theta & -\sin \theta \\
	\sin \theta & \cos \theta \\
\end{array}
\right ]
\left [
\begin{array}{cc}
 -\sin \beta \\
 \cos \beta  \
\end{array} \right ]=
\left [
\begin{array}{cc}
 -\sin (\theta + \beta) \\
 \cos (\theta + \beta)  \
\end{array} \right ].
\]

The condition (\ref{eq2}) $ {\cal R} \vec{t}\cdot \vec{j}=0$ implies that $\cos (\theta + \beta)=0$, that is, 
\[
\theta + \beta =(2n -1)\fracd{\pi}{2},~n=1,2. 
\]
So, we have two different solutions:
\[
\begin{array}{c}
\beta_s = \fracd{\pi}{2} - \theta, \\
\beta_i = \fracd{3\pi}{2} - \theta.
\end{array}
\]
Suppose that the disc with radius $a$ rotates with constant angular velocity $\theta=\omega t$. Therefore it is possible to write:

\[
\begin{array}{c}
\beta_s = \fracd{\pi}{2} - \omega t, \\
\beta_i = \fracd{3\pi}{2} - \omega t.
\end{array}
\]

Now, the general expression of the position vector on $OXY$ is

\[\vec{R}=a
\left [
\begin{array}{cc}
	\cos \theta & -\sin \theta \\
	\sin \theta & \cos \theta \\
\end{array}
\right ]
\left [
\begin{array}{cc}
 \cos \beta \\
 \sin \beta  \
\end{array} \right ]=a
\left [
\begin{array}{cc}
 \cos (\theta + \beta) \\
 \sin (\theta + \beta)  \
\end{array} \right ].
\]

It follows that
\[
Y_s = \vec{R}(\beta_s) \cdot \vec{j} = a \sin (\theta + \beta_s ).
\]

Since $\beta_s = \fracd{\pi}{2} - \theta$ it comes:
\[
Y_s = a \sin \left ( \fracd{\pi}{2} \right ) =a.
\]
Similarly, $Y_i=-a$ and $H=Y_s - Y_i = 2a$.

\esp

A more interesting case arises when the disc rotates around a point that does not coincides with the center. Take the case of a disc rotating around a fixed point on the contour:
\[
\vec{r_0}=\left [
\begin{array}{c}
1 \\
0 
\end{array}\right ],~\theta (t)=\omega t.
\]
\begin{figure}[!ht]
\centering
\includegraphics[scale=0.5]{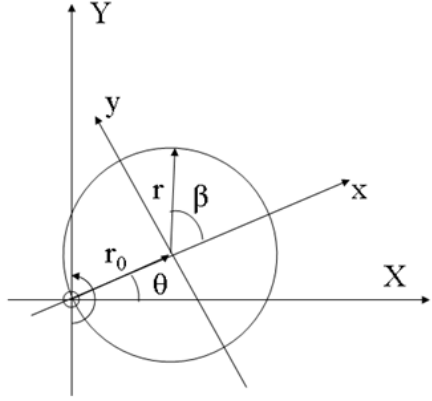}
\caption{A circle rotating around a point of the contour}
\end{figure}

There is no loss of generality if we take the $x$ axis fixed on the disc through the rotation center. Following the reasoning exposed above, we have successively:

\[ 
\vec{r}= a \left [
\begin{array}{cc}
 \cos \beta \\
 \sin \beta  \
\end{array} \right ],
\vec{t}= \left [
\begin{array}{cc}
 -\sin \beta \\
 \cos \beta  \
\end{array} \right ],
\]

\[{\cal R} \vec{t}=
\left [
\begin{array}{cc}
 -\sin (\theta + \beta) \\
 \cos (\theta + \beta)  \
\end{array} \right ].
\]

The condition (\ref{eq2}) $ {\cal R} \vec{t}\cdot \vec{j}=0$ implies that 
\[
\begin{array}{c}
\beta_s = \fracd{\pi}{2} - \theta, \\
\beta_i = \fracd{3\pi}{2} - \theta.
\end{array}
\]

The general expression of the position vector on $OXY$ is

\[\vec{R}={\cal R}(\vec{r_0} + \vec{r})
=\left [
\begin{array}{cc}
	\cos \theta & -\sin \theta \\
	\sin \theta & \cos \theta \\
\end{array}
\right ]
\left (
a \left [
\begin{array}{c}
1 \\
0
\end{array}
\right ] + a \left [
\begin{array}{c}
 \cos \beta \\
 \sin \beta  
\end{array} \right ]
\right ).
\]

Its projection on the $Y$ axis is
\[
\vec{R} \cdot \vec{j} = a \left [ \sin \theta + \sin (\theta + \beta) \right ] =Y.
\]

Substituting $\beta$ by $\beta_s$ and $\beta_i$ we find

\[
\begin{array}{c}
Y_s = Y(\beta_s)=a[\sin (\theta) + 1], \\
Y_i = Y(\beta_i)=a[\sin (\theta) - 1]
\end{array}
\]

and $H=Y_s - Y_i = 2a$.

The printed image over the film that winds down in front of the slot with velocity $v=v_0=$constant is obtained from the system

\[ 
\begin{array}{cl}
\fracd{dY_s}{dt}=& a\cos \theta \frac{d\theta}{dt}=a\omega \cos (\omega t),\\
\fracd{dz}{dt} =& v_0 .	
\end{array}
\]

Using the conditions $Y_s (0)=a$ and $z(0)=0$ the solutions of the equations above are:
\[ 
\begin{array}{cl}
Y_s = & a \left ( 1 + \sin (\omega t) \right ),\\
z =& v_0 t. 	
\end{array}
\]
Solving this system, by elimination of the variable $t$ we obtain:
\[ 
Y_s= a\left ( 1 + \sin \left( \fracd{\omega}{v_0} z \right) \right ).
\]

Similarly we also obtain:

\[
Y_i=Y_s - H= a \left ( -1 + \sin \left ( \fracd{\omega}{v_0} z  \right ) \right ). 	
\]

The kinematic image has the shape of a wave (see figure 11).
\begin{figure}[!ht]
\centering
\includegraphics[scale=0.5]{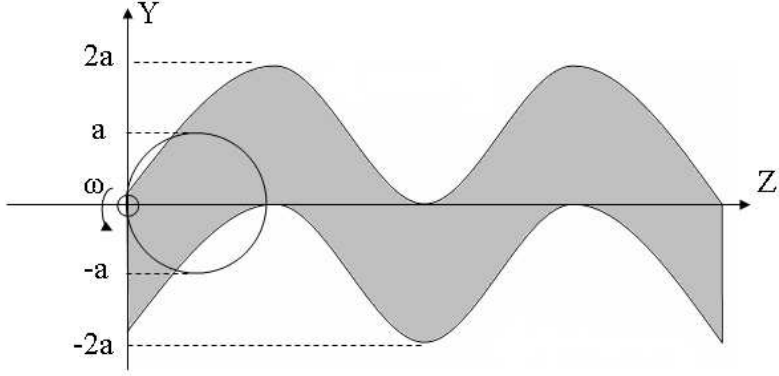}
\caption{Kinematic image of a circle rotating around a point on the contour}
\end{figure}

\renewcommand{\theequation}{5.\arabic{equation}}
\setcounter{equation}{0}

{\bf 5. THE CASE OF THE ELIPSE.} Consider the position vector of an arbitrary point on the contour of an elipse centered in the origin of a system $OXY$ with semi-axis $a > b$:

\[ \vec{r}= \left [
\begin{array}{cc}
a \cos \beta \\
b \sin \beta  \
\end{array} \right ].
\]

Its tangent vector is given by:

\[ \vec{t}= \left [
\begin{array}{cc}
-a \sin \beta \\
b \cos \beta  
\end{array} \right ].
\]

The rotation of tangent vector by an angle $\theta$ is:

\[ {\cal R} \vec{t}= \left [
\begin{array}{cc}
-a \cos \theta \sin \beta -b \sin \theta \cos \beta \\
-a \sin \theta \sin \beta + b \cos \theta \cos \beta  
\end{array} \right ].
\]

The tangency condition (\ref{eq2}) leads to

\begin{equation} \label{5.1}
a \sin \theta \sin \beta = b \cos \theta \cos \beta.
\end{equation}

The general expression of the position vector referred to $OXY$ is

\[\vec{R}=
\left [
\begin{array}{cc}
	\cos \theta & -\sin \theta \\
	\sin \theta & \cos \theta \\
\end{array}
\right ]
\left [
\begin{array}{c}
a \cos \beta \\
b \sin \beta
\end{array}
\right ]
=\left [
\begin{array}{c}
a \cos \theta \cos \beta - b \sin \theta \sin \beta \\
a \sin \theta \cos \beta + b \cos \theta \sin \beta
\end{array}
\right ]
\]
and its projection on the axis $OY$ is
\begin{equation}\label{5.2}
Y=\vec{R}\cdot \vec{j}=a \sin \theta \cos \beta + b \cos \theta \sin \beta.
\end{equation}

Rather than try to write $\beta$ as a function of $\theta$ let explore useful relationships between (\ref{5.1}) and (\ref{5.2}). In fact, squaring (\ref{5.2}) we have
\[
Y^2 = a^2 \sin^2 \theta \cos^2 \beta + b^2 \cos^2 \theta \sin^2 \beta + 2ab \sin \theta \cos \beta \cos \theta \sin \beta.
\]
Using the tangency condition (\ref{5.1}) we obtain successively
\[
\begin{array}{ll}
Y^2_{s,i} &= a^2 \sin^2 \theta \cos^2 \beta + b^2 \cos^2 \theta \sin^2 \beta + 2a^2  \sin^2 \theta \sin^2 \beta \\
& = a^2 \sin^2 \theta + b^2 \cos^2 \theta \sin^2 \beta + a^2  \sin^2 \theta \sin^2 \beta,
\end{array}
\]
that is,
\begin{equation} \label{5.3}
Y^2_{s,i} = a^2 \sin^2 \theta + (a^2  \sin^2 \theta + b^2 \cos^2 \theta)\sin^2 \beta.
\end{equation}
Similarly we get: 
\begin{equation}\label{5.4}
Y^2_{s,i} = b^2 \cos^2 \theta + (a^2  \sin^2 \theta + b^2 \cos^2 \theta)\cos^2 \beta.
\end{equation}
Adding (\ref{5.3}) and (\ref{5.4}) we conclude that
\[
Y^2_{s,i} = a^2  \sin^2 \theta + b^2 \cos^2 \theta,
\]
which means that
\begin{equation}
\begin{array}{ll}
Y_s = & \sqrt{a^2 \sin^2 \theta + b^2 \cos^2 \theta},\\
Y_i = & -\sqrt{a^2 \sin^2 \theta + b^2 \cos^2 \theta}.
\end{array}
\end{equation}

Since $\theta=\omega t$ and $z=v_0 t$, the printed image on the film is obtained from the equations:
\begin{equation}
\begin{array}{ll}
Y_s = & \sqrt{a^2 \sin^2\left ( \displaystyle{\frac{\omega}{v_0}z}\right) + b^2 \cos^2 \left (\displaystyle{\frac{\omega}{v_0}z}\right ) },\\
 & \\
Y_i = & -\sqrt{a^2 \sin^2 \left ( \displaystyle{\frac{\omega}{v_0}z} \right ) + b^2 \cos^2 \left ( \displaystyle{\frac{\omega}{v_0}z} \right )}.
\end{array}
\end{equation}
The figure 12 shows the printed image on a moving film of an elipse rotating around its center with constant velocity $v_0$:
\begin{figure}[!ht]
\centering
\includegraphics[scale=0.5]{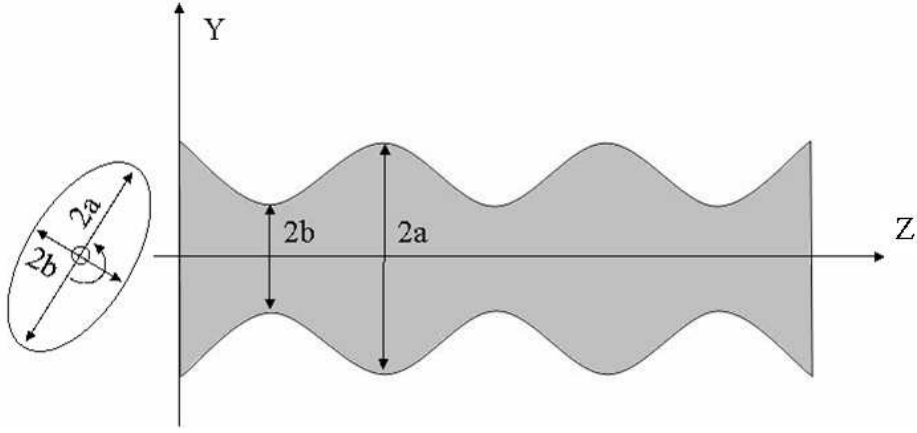}
\caption{Kinematic image of an elipse rotating around its center}
\end{figure}

\esp

\renewcommand{\theequation}{6.\arabic{equation}}
\setcounter{equation}{0}

{\bf 7. THE CASE OF REGULAR POLYGONS.} The techniques presented in the section 3 cannot be applied to this case since we are dealing with plane figures with corners, non-smooth contour. We will see that for this case the profile projected on the moving screen is defined by the vertices or by the corresponding position vectors. 

For sake of simplicity let us call  $\vec{P}$ the vector $\overrightarrow{OP}$  where $O$ is the origin of the coordinate system and $P$ a general vertex of the polygon. 
Now given a vector $\vec{P}={[x,y]}^T$ , the $Y$-component after a rotation $\theta$ is given by:
 
\[
f_{\vec{P}} (\theta)=x \sin \theta + y \cos \theta.
\]

Let us call this function {\it RotProj}. Consider the polygon inscribed in a circle whose center coincides with the origin of the $X-Y$ axis system. Assume that the polygon rotates around the point coinciding with the center of the circle with constant angular velocity $\omega$ and the film moves with constant velocity $v_0$. 

{\bf Case 1: the square.} Let us consider a square with side equals $a$ as shown in the figure 13.
\begin{figure}[!ht]
\centering
\includegraphics[scale=0.5]{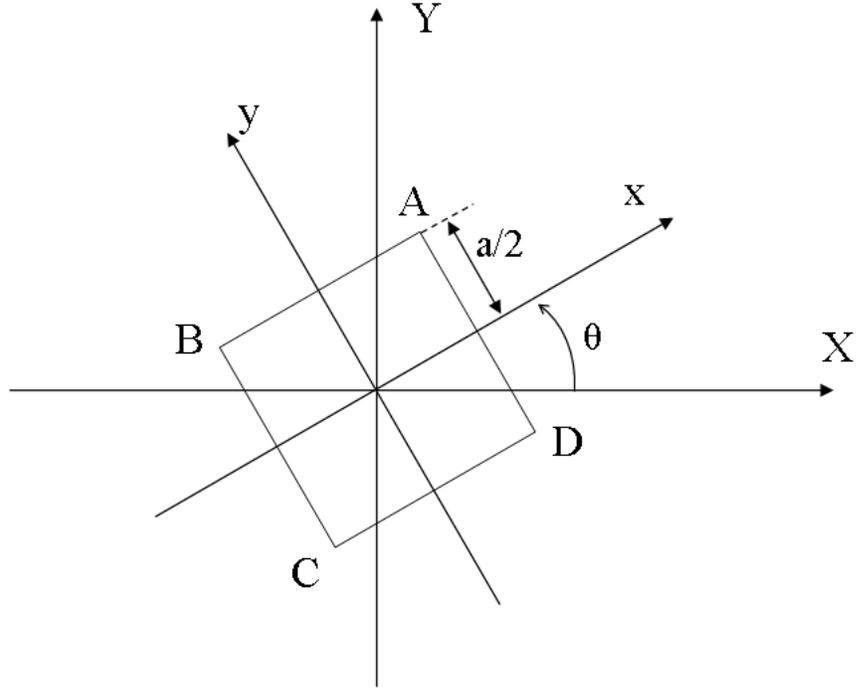}
\caption{Image of a square}
\end{figure}

For $\theta=0$ the vertices are given by:

\[
\begin{array}{l}
A=\frac{a}{2} (1,1), \\
B=\frac{a}{2} (-1,1) ,\\
C=\frac{a}{2} (-1,-1),\\
D=\frac{a}{2} (1,-1).
\end{array}
\]

The corresponding functions RotProj read:

\[
\begin{array}{l}
f_{\vec{A}}(\theta)=\frac{a}{2} (\sin \theta + \cos \theta),\\
f_{\vec{B}}(\theta)=\frac{a}{2} (\sin \theta - \cos \theta), \\
f_{\vec{C}}(\theta)=\frac{a}{2} (-\sin \theta - \cos \theta), \\
f_{\vec{D}}(\theta)=\frac{a}{2} (-\sin \theta + \cos \theta).
\end{array}
\]

Now when the square rotates of an angle $\theta$ within the interval $0\leq \theta < \pi /2$, the vertices defining the upper and bottom contours of the profile projected on the film are $A$   and  $C$ respectively. The vertices $B$ and $D$ in that case are "hidden" (see figure 14).
\begin{figure}[!ht]
\centering
\includegraphics[scale=0.5]{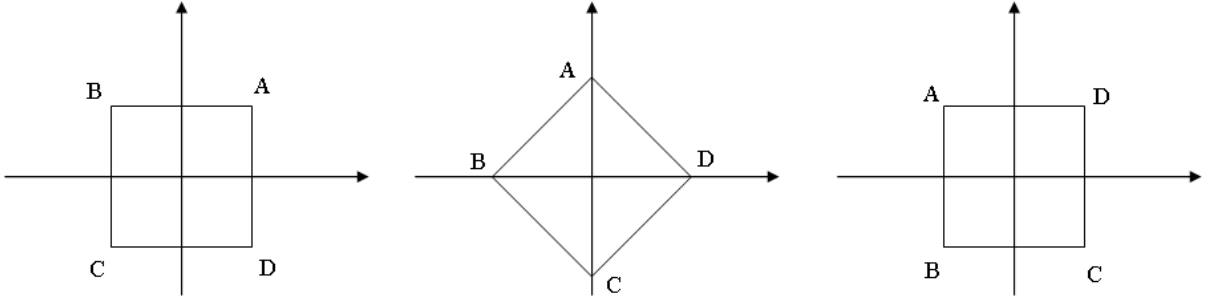}
\caption{Successive positions of a rotating square}
\end{figure}
The expressions for the upper and bottom contours are given by  $Y_s = f_{\vec{A}}(\theta)=a/2 (\sin \theta + \cos \theta )$ and $Y_i = -Y_s$  respectively. 

For the case where $\theta$ belongs to the interval $\pi /2 \leq \theta < \pi$, the control vertices are $D$ and $B$ while $A$ and $C$ remain hidden.

Following this procedure till the rotation completes a complete cycle $(2\pi)$ the expression from the upper and bottom contours are given by:

\[
Y_s = \cases{ \frac{a}{2} (\sin \theta + \cos \theta), &  if $ 0 \leq \theta < \pi /2$, \cr 
             \frac{a}{2} (\sin \theta - \cos \theta),  &  if $\pi /2 \leq \theta < \pi$, \cr 
             \frac{a}{2} (-\sin \theta - \cos \theta), & if  $\pi \leq \theta < 3\pi /2$, \cr
             \frac{a}{2} (-\sin \theta + \cos \theta), & if  $3 \pi /2 \leq \theta < 2\pi$, \cr} 
\]
and  $Y_i = - Y_s$.   

The image projected on the film, that is, the kinematic image, is depicted in the figure 15.
\begin{figure}[!ht]
\centering
\includegraphics[scale=0.4]{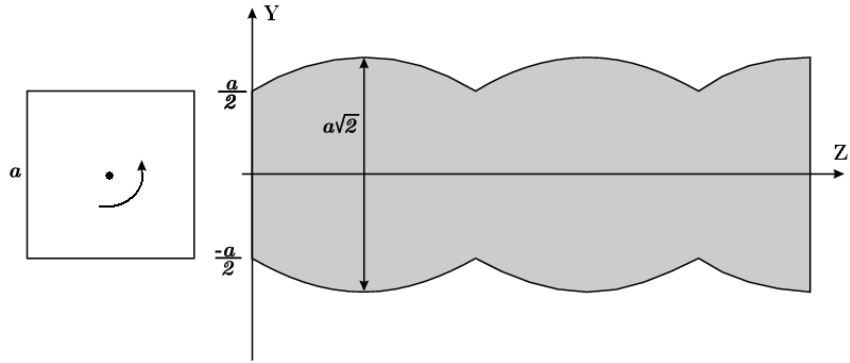}
\caption{Kinematic image of a square  rotating around its center with constant angular velocity. Film moving with constant speed $v=v_0$}
\end{figure}

\newpage

{\bf Case 2: The Equilateral Triangle.} Consider now an equilateral triangle with side equal $a$ as shown in the figure 16.
\begin{figure}[!ht]
\centering
\includegraphics[scale=0.5]{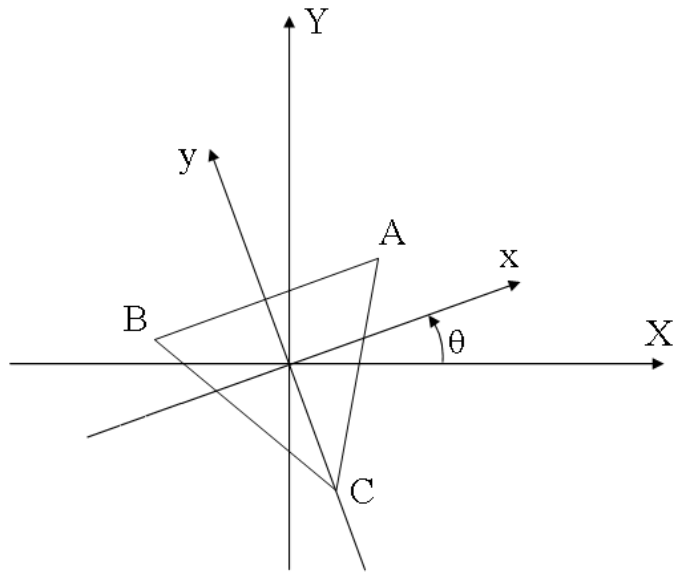}
\caption{Image of an equilateral triangle}
\end{figure}

For $\theta=0$ the vertex coordinates are given by:

\[
\begin{array}{lll}
A=\displaystyle{\frac{a\sqrt{3}}{6}} (\sqrt{3},1), & B=\displaystyle{a\frac{\sqrt{3}}{6}} (-\sqrt{3},1), & C=\displaystyle{\frac{a\sqrt{3}}{3}}(0,-1).
\end{array}
\]

The corresponding projections are given by the functions:
\[
\begin{array}{l}
f_{\vec{A}}(\theta)=\displaystyle{\frac{a\sqrt{3}}{6}} (\sqrt 3 \sin \theta + \cos \theta ), \\
f_{\vec{B}}(\theta)=\displaystyle{\frac{a\sqrt{3}}{6}} (-\sqrt 3 \sin \theta + \cos \theta ), \\
f_{\vec{C}}(\theta)=-\displaystyle{\frac{a\sqrt{3}}{3}} \cos \theta.
\end{array}
\]

When the triangle rotates of an angle $\theta$ comprised within the interval $0\leq \theta \leq \pi /3$, the upper contour is determined by $\vec{A}$   and the bottom contour by $\vec{C}$ . In that case the vertex $B$ is hidden (see figure 17).
\begin{figure}[!ht]
\centering
\includegraphics[scale=0.5]{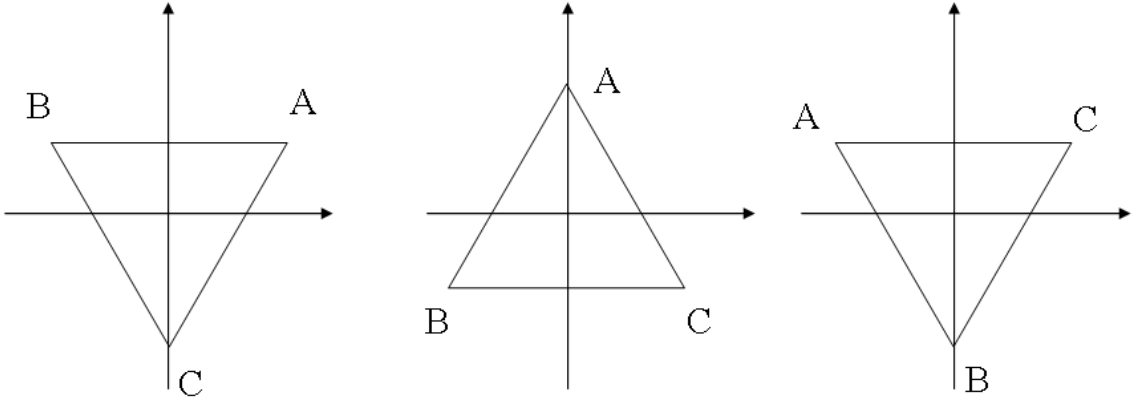}
\caption{Sucessive images of an equilateral triangle rotanting around its center}
\end{figure}

Therefore, within this angular interval the upper and bottom projections are given by 
$Y_s = f_{\vec{A}}(\theta)=\displaystyle{\frac{a\sqrt{3}}{6}} (\sqrt 3 \sin \theta + \cos \theta )$
and  $Y_i = f_{\vec{C}}(\theta)=-\displaystyle{\frac{a\sqrt{3}}{3}} \cos \theta $, respectively.

For the next interval $\pi /3 \leq \theta \leq 2 \pi /3$ the projection is controlled by $A$ and $B$ with $C$ hidden.

The table below shows the intervals composing a complete cycle $(2\pi)$, the corresponding control vertices and the hidden vertex as well. 

\begin{tabular}{|c|c|c|c|}
\hline
Angular interval & Upper control vertex  & Bottom control vertex  & Hidden vertex\\ \hline
\hline \hline
$0 \leq \theta \leq \pi/3 $      & $A$ & $C$   & $B$        \\ \hline 
$\pi/3 \leq \theta \leq 2\pi/3$  & $A$ & $B$   & $C$          \\ \hline 
$2\pi/3 \leq \theta \leq \pi$    & $C$ & $B$   & $A$   \\  \hline  
$\pi \leq \theta \leq 4\pi/3$    & $C$ & $A$   & $B$        \\ \hline
$4\pi/3 \leq \theta \leq 5\pi/3$ & $B$ & $A$   & $C$          \\ \hline
$5\pi/3 \leq \theta \leq 2\pi$   & $B$ & $C$   & $A$ \\   \hline 
\end{tabular}

\esp

Using the table above the respective functions {\it RotProj} give the upper and the bottom contour curves. They are:

\[
Y_s = \cases{ 
\displaystyle{\frac{a\sqrt{3}}{6}} (\sqrt 3 \sin \theta + \cos \theta ), &  if $0 \leq \theta \leq 2\pi/3 $,  \cr 
-\displaystyle{\frac{a\sqrt{3}}{3}} \cos \theta ,  &  if  $2\pi/3 \leq \theta \leq 4\pi/3$,\cr 
 \displaystyle{\frac{a\sqrt{3}}{6}} (-\sqrt 3 \sin \theta + \cos \theta )   , & if  $4\pi/3 \leq \theta \leq 2\pi$, \cr} 
\]
and
\[
Y_i = \cases{ -\displaystyle{\frac{a\sqrt{3}}{3}} \cos \theta , & if $0 \leq \theta \leq \pi/3 $,  \cr
\displaystyle{\frac{a\sqrt{3}}{6}} (-\sqrt 3 \sin \theta + \cos \theta ), &  if $\pi/3 \leq \theta \leq \pi $,  \cr 
\displaystyle{\frac{a\sqrt{3}}{6}} (\sqrt 3 \sin \theta + \cos \theta )  , & if  $\pi \leq \theta \leq 5\pi/3$, \cr
-\displaystyle{\frac{a\sqrt{3}}{3}} \cos \theta ,  &  if  $5\pi/3 \leq \theta \leq 2\pi$.\cr 
 } 
\]

The image printed on the film is given by figure 18.
\begin{figure}[!ht]
\centering
\includegraphics[scale=0.5]{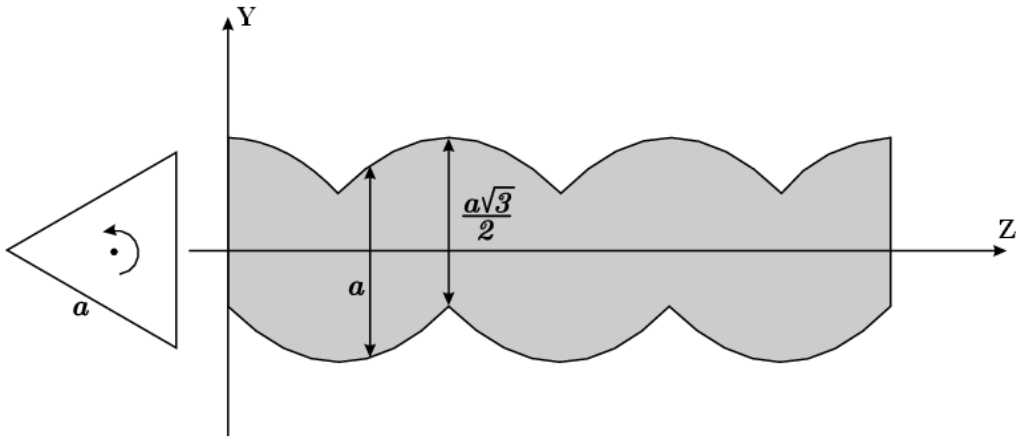}
\caption{Kinematic image of an equilateral triangle rotanting around its center}
\end{figure}

Exercice: Find the shape of the kinematic image on the film of a regular polygon having $n$ sides. What is the difference between the image of the polygons consisting of even and odd sizes?

\newpage

\renewcommand{\theequation}{7.\arabic{equation}}
\setcounter{equation}{0}

{\bf 7. THE INVERSE PROBLEM.} We have discussed the solution of the direct problem. That is, find the image of a bi-dimensional rotating solid projected on a moving film. The general inverse problem is more complex and has yet to be solved for the general case. However, some conclusions taken from the direct problem can help the formulation of the inverse problem. An example is the case of regular polygons which was discussed in the previous section. 

Indeed, consider an image generated by the rotation of a regular polygon with constant angular velocity. Suppose also that the film moves with a constant speed. We want to know which polygon has generated that image. To solve this problem, suppose that the polygon rotates around its center of mass. If this centre is located at the origin of a system of Cartesian coordinates, then the shadow will be limited above and below by periodic curves. We can ensure, in advance, that if those curves are symmetrical with respect to the $x$ axis then the polygon has a number of even sizes. If the curves are shifted, the polygon has an odd number of sides. To know the exact number of sides of the polygon it is enough to know the relationship between the values $m$, meaning the minimum distance from the superior curve to the $x$ axis, and $M$, meaning the maximum distance from the superior curve to the $x$ axis. The value of $m$ is given by  polygon's apotema and $M$ is the distance of the center to the vertex of the polygon (see figure 19).
\begin{figure}[!ht]
\centering
\includegraphics[scale=0.5]{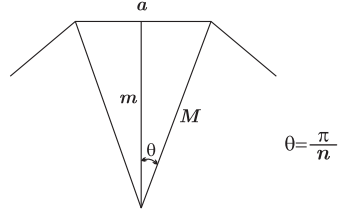}
\caption{Image of a regular polygon}
\end{figure}

Clearlyt $\theta=\displaystyle{\frac{\pi}{n}=\arccos{\left(\frac{m}{M}\right)}}$, where $n$ is the number of edges of the regular polygon. Thus

\[
n=\frac{\pi}{\arccos{\left(\frac{m}{M}\right)}}.
\]

Por instance, if we do not know in advance which the regular polygon generated the figure 15, by inspection of the image we may conclude immediately the regular polygon has a number of even sizes due the symmetry ot the upper and lower bounds. To find out the exact number of edges it is enough to calculate
 $m$ and $M$. In fact, $m=\displaystyle{\frac{a}{2}}$ and $M=\displaystyle{ \frac{a\sqrt 2 }{2} }$. Therefore,
\[
n=\frac{\pi}{\arccos{\left(\frac{a/2}{a\sqrt 2 /2}\right)}}=\frac{\pi}{\pi /4}=4,
\]
that is, the polygon is a square.

For the case of figure 19, $m=\displaystyle{\frac{a \sqrt 3}{6}}$ and $M=\displaystyle{\frac{a \sqrt 3}{3}}$. We conclude that $n=3$ and the polygon is a equilateral triangle.

Note that the smaller the difference $(M-m)$, the greater the number of sides of the regular polygon. If $m \rightarrow M$ then
\[
\lim_{m \rightarrow M}  \frac{\pi}{ \arccos{ \left(\frac{m}{M}\right) } } =+\infty,
\] 
that is, in the limit, we have a rotating circle.

\newpage

\renewcommand{\theequation}{8.\arabic{equation}}
\setcounter{equation}{0}

{\bf 8. FINAL REMARKS.} We have discussed in this short paper a simple problem. We believe however that it clarifies and illustrates some basic questions arising in direct and inverse problems. It becomes very clear, for instance, why convexity of domains of definition are so important for the uniqueness of the solutions, particularly for the case of inverse problems.  

The kinematic image can also be used as a tool to stimulate the abstract representation of the shape of two-dimensional objects. In that sense it could be an interesting tool for educational purposes. 
Also more complex motion of the film and of the rotating object can be introduced. It is possible that the identification of intricate contours can take advantage of the comparison between two or more kinematic images corresponding to different motions. 
We believe that this problem is a starting point for several other setups including three-dimensional objects.

\end{document}